\documentclass[twoside,12pt]{article}
\usepackage{amstex,amssymb,a4wide,doublespace,epic,eepic}
\usepackage{babel}
\newtheorem{theorem}{Theorem}[section]
\newtheorem{lemma}[theorem]{Lemma}
\newtheorem{proposition}[theorem]{Proposition}

\newtheorem{definition}[theorem]{Definition}
\newtheorem{corollary}[theorem]{Corollary}
\newtheorem{example}[theorem]{Example}

\newtheorem{remark}[theorem]{Remark}
\newcommand{\bproof}{\noindent{\bf Proof: }}
\newcommand{\eproof}{\hfill $\Box$\\}

\newcommand{\Span}{\mbox{\rm span}}
\newcommand{\dist}{\mbox{\rm dist}}
\newcommand{\GL}{\mbox{\rm GL}}
\newcommand{\U}{\textrm{U}}
\newcommand{\tN}{\textrm{N}}
\newcommand{\red}{\textrm{red}}
\newcommand{\sr}{\textrm{sr}}

\newcommand{\w}{\omega}
\newcommand{\e}{\varepsilon}
\newcommand{\bC}{{\mathbb C}}
\newcommand{\bR}{{\mathbb R}}
\newcommand{\bT}{{\mathbb T}}
\newcommand{\bZ}{{\mathbb Z}}
\newcommand{\Bhat}{\hat{B}}
\newcommand{\bN}{{\mathbb N}}
\newcommand{\frA}{\mathfrak{A}}

\renewcommand{\ell}{l}
\renewcommand{\Re}{\mbox{\rm Re}}
\renewcommand{\Im}{\mbox{\rm Im}}
 
\begin{document}
\markboth{Ken Dykema, Uffe Haagerup, Mikael R\o rdam}{The Stable Rank of
Some Free Products}
 
\title{The Stable Rank of Some Free Product $C^*$-Algebras}
 
\author{Ken Dykema \and Uffe Haagerup \and Mikael R\o rdam\\Dept.\ of
Mathematics and Computer Science,\\Odense University,
DK-5230 Odense M, Denmark}

\maketitle

\begin{abstract}
It is proved that the reduced group $C^*$-algebra $C^*_{\red}(G)$ has
stable rank one (i.e.\ its group of invertible elements is a dense
subset) if $G$ is a discrete group arising as a free product $G_1*G_2$
where $|G_1|\ge 2$ and $|G_2|\ge 3$. This follows from a more general
result where it is proved that if $(\frA,\tau)$ is the reduced free
product of a family $(A_i,\tau_i)$, $i\in I$, of unital $C^*$-algebras
$A_i$ with normalized faithful traces
$\tau_i$, and if the family satisfies the
Avitzour condition (i.e.\ the traces, $\tau_i$, are not too lumpy in a
specific sense), then $\frA$ has stable rank one. 
\end{abstract}

\section*{Introduction}

It is an open problem if every finite, simple $C^*$-algebra has stable
rank one. Recall that a unital $C^*$-algebra $A$ is said to have stable
rank one if the group of invertible elements in $A$ is a norm dense
subset of $A$.

The notion of stable rank was introduced by M.~Rieffel in \cite{Rf} with
the purpose of establishing what one might call non-stable $K$-theory
results for certain concrete $C^*$-algebras, most notably the irrational
rotation $C^*$-algebras. On a more speculative note, stable rank (which
associates a number in $\{1,2,3,\dots\}\cup\{\infty\}$ to every
$C^*$-algebra) should measure the (non-commutative) dimension of the
$C^*$-algebra, with stable rank equal to one corresponding to dimension
0 or 1. (It has later turned out that different definitions of
dimensions, that agree in the ``commutative'' case, generalize to
dimension concepts for $C^*$-algebras which do not agree.)

Some of the non-stable $K$-theory results, obtained in \cite{Rf} and
\cite{Rf1}, for
unital $C^*$-algebras $A$ of stable rank one are as follows.
The three relations on projections in $A$ (or in a matrix algebra over
$A$), Murray--von Neumann equivalence, unitary equivalence and homotopy
equivalence, are the same. Moreover, if $p,q$ are projections in $A$ (or
in a matrix algebra over $A$), such that $[p]_0=[q]_0$ in $K_0(A)$, then
$p$ and $q$ are equivalent with respect to either of the three relations
mentioned above. Also, the natural group homomorphism
\[
\U(A)/\U_0(A)\to K_1(A),
\]
where $\U(A)$ is the group of unitary elements in $A$ and $\U_0(A)$ its
connected component containing the unit of $A$, is an isomorphism. 

Another property of $C^*$-algebras of stable rank one can be found in
\cite{FR}. It is observed in that paper that if $A$ is a $C^*$-algebra
of stable rank one, then each normal element in $A$ can be approximated
by normal elements in $A$ with 1-dimensional spectrum. This again is
used to prove that there exists a function $f: \bR^+\to \bR^+$
(independent of $A$) which is continuous at 0 and with $f(0)=0$, such
that
\[
\dist\big(a,\tN(A)\big) \le f\big(\|a^*a-aa^*\|\big)
\]
for every $a\in A$ with $\|a\|\le 1$. (Here $\tN(A)$ denotes the set of
normal elements in $A$.)

The main result of this paper (Theorem \ref{thm-3.8}) states that the
reduced free product of any pair of unital $C^*$-algebras $A_1$ and
$A_2$ with faithful normalized traces $\tau_1$ and $\tau_2$ has stable
rank one if there exist unitary elements $u\in A_1$, $v,w\in A_2$ such
that $0=\tau_1(u)=\tau_2(v)=\tau_2(w)=\tau_2(w^*v)$ (the Avitzour
condition). This result applies in particular to the reduced group
$C^*$-algebras $C^*_{\red}(G)$ when $G=G_1*G_2$ for some groups $G_1$
and $G_2$ satisfying $|G_1|\ge 2$ and $|G_2|\ge 3$ (see Corollary
\ref{cor-3.9}). It follows in particular that $C^*_{\red}(F_n)$, $2\le
n\le\infty$, and the Choi algebras $C^*_{\red}(\bZ_n*\bZ_m)$,
$n\ge 2$ and $m\ge 3$, have stable rank one.

Our result answers a question of Marc Rieffel~\cite{RieffelZZKthCrPr}, showing
that every
projective module over $C^*_{\red}(F_n)$, $2\le n\le \infty$, is
free. Indeed, if $p$ is a projection in a matrix algebra over
$C^*_{\red}(F_n)$, then by~\cite{PV} we have
$[p]_0 = k\cdot [1]_0$ in $K_0(C^*_{\red}(F_n))$
for some $k\in\bN$. Now, if $q$ is the free module over
$C^*_{\red}(F_n)$ of dimension $k$ (viewed as a projection in the same
matrix algebra over $C^*_{\red}(F_n)$ as $p$ belongs to), then
$[p]_0=[q]_0$. Since $C^*_{\red}(F_n)$ has stable rank one, it follows
that $p$ and $q$ are equivalent, and hence $p$ is free.
For comparison, it has been known for some time that the group ring, $\bC F_n$,
is a (left and right) free
ideal ring, (fir), see~\cite{CohnZZFreeRings}, and hence every submodule
of a free module over $\bC F_n$ is free.
That $\bC F_n$ is a fir follows because, from \cite{Bergman} Corollary~2.12,
the free product (also called coproduct) of firs is a fir, and because
$\bC\bZ$ is a principal ideal domain, thus a fir.

It also follows (see~\cite{PV}) that
$\U(C^*_{\red}(F_n))/\U_0(C^*_{\red}(F_n))$
is naturally isomorphic to $\bZ^n$.
In particular, with $\lambda:F_n\to C^*_{\red}(F_n)$ the left regular
representation, $\lambda(g)$ is connected in $\U(C^*_{\red}(F_n))$ to $1$ if
and
only if $g$ belongs to the commutator subgroup of $F_n$.

The strategy of the proof of Theorem \ref{thm-3.8} follows in parts  (Lemmas
\ref{lemma-3.3} and \ref{lemma-3.4}) the work of the second named author in
\cite{Ha}.
A crucial ingredient of the proof is the result of the third named
author~\cite{Ro}, that if $\frA$ is a C$^*$--algebra whose group of invertibles
is not dense in $\frA$, then there is an element, $b\in\frA$ of norm $1$ whose
distance to the invertibles is equal to $1$.

In order to emphasize the main ideas of the proof, we will first, in Section
\ref{sec-new2}, go through the proof that the invertibles are dense in
$C^*_{\red}(F_2)$.
In Section \ref{sec-2} we give some preliminaries for the proof of the more
general result that the reduced free product of two unital C$^*$--algebras
satisfying the Avitzour conditions has stable rank one, and in Section
\ref{sec-3} we prove this theorem.
In Section \ref{sec-4} we derive some conditions under which the
Avitzour condition is satisfied (see Proposition \ref{prop-4.1}). In
particular, the Avitzour condition is satisfied if $A_1$ and $A_2$
both contain unital abelian subalgebras which are non-atomic with
respect to the traces $\tau_1$ and $\tau_2$. Section \ref{sec-5}
contains a brief discussion of the structure of more general reduced
free products.

\section{The proof of a special case}
\label{sec-new2}

\begin{theorem}
\label{thm-new2.1}
The reduced group C$^*$--algebra $C^*_{\red}(F_2)$ has stable rank one.
\end{theorem}

\bproof
Write $F_2=\langle a,b\rangle$, i.e.\ $F_2$ is freely generated by $a$ and $b$,
and write $\frA=C^*_{\red}(F_2)$.
Thus $\frA$ is generated by the set of left translators,
$\{\lambda_g\mid g\in F_2\}$.
Let $\tau$ denote the canonical, faithful, tracial state
on $\frA$.
We then have the inner product $\langle w,z\rangle=\tau(z^*w)$ and we denote
$\|z\|_2=\langle z,z\rangle^{1/2}$.
Note that $(\lambda_g)_{g\in F_2}$ is an orthonormal basis for (the closure of)
$\frA$ with respect to this inner product.

Suppose for contradiction that $\frA$ has stable rank strictly
greater that $1$.
Then by Theorem 2.6 of \cite{Ro}, there must be $x\in\frA$ such that
$\|x\|=1$ and the distance in norm from $x$ to the invertibles of
$\frA$, (denoted $\GL(\frA)$), is $1$.
But then $\|x\|_2<1$ because $\|x\|_2=1$ would imply that $x$
be unitary.
Hence we can find a finite linear combination of left translators,
$y=\sum_{j=1}^n\alpha_j\lambda_{g_j}$, such that $\|y\|_2$ is strictly less
than the distance from $y$ to $\GL(\frA)$.
Taking $k\in\bN$ and considering each $b^kg_jb^{-k}$ as a reduced word in $a$,
$b$ and their inverses, we see that
there is $k$ large enough so that
for every $j$, $b^kg_jb^{-k}$ when reduced begins and
ends with $b$ or $b^{-1}$.
Thus we see that there is no cancellation when we multiply
$(ab^kg_jb^{-k}a)(ab^kg_{j'}b^{-k}a)$, for any $j$ and $j'$.
This shows that when $u=\lambda_{ab^k}$ and $v=\lambda_{b^{-k}a}$ we have
\begin{equation}
\label{2normpowers}
\|(uyv)^m\|_2=\|uyv\|_2^m,\qquad m\in\bN.
\end{equation}

In \cite{Ha} it was shown that for finite linear combinations
of left translators from $F_2$, the operator norm can be bounded in terms of
the two--norm.
Indeed from Lemma~1.5 of that paper it follows that if
$z=\sum_{j=1}^p\beta_j\lambda_{h_j}$ and if $N$ is the maximum of the lengths
of the words $h_j$, (as reduced words in $a$, $b$ and their inverses), then
\[
\|z\|\le2(N+1)^2\|z\|_2.
\]
Now we can apply the estimates from (\ref{2normpowers}) to obtain an upper
bound
on the spectral radius of $uyv$, namely, letting $N$ be the maximum of the
lengths of the words $ab^kg_1b^{-k}a,\ldots,ab^kg_nb^{-k}a$, we have
\begin{eqnarray*}
r(uyv) &=& \limsup_{m\to\infty}\|(uyv)^m\|^{1/m} \\
&\le& \limsup_{m\to\infty}(2(mN+1)^2\|uyv\|_2^m)^{1/m}=\|uyv\|_2=\|y\|_2.
\end{eqnarray*}
But the distance from $uyv$ to $\GL(\frA)$ is clearly no greater than the
spectral radius of $uyv$ and the distance from $y$ to $\GL(\frA)$ is equal to
the distance from $uyv$ to $\GL(\frA)$, so the inequality
$r(uyv)\le\|y\|_2$
gives a contradiction to the choice of $y$.
\eproof

\section{Preliminaries for the general case}
\label{sec-2}

\subsection{Standard orthonormal basis}
\label{subsec-2.1}
Let $A$ be a unital $C^*$-algebra with a faithful normalized trace
$\tau$. Consider the corresponding Euclidean structure:
\begin{eqnarray*}
\langle a,b\rangle &=& \tau(b^*a),\qquad a,b\in A,\\
\|a\|_2 &=& \langle a,a\rangle^{\frac12},\qquad a\in A.
\end{eqnarray*}

A subset $X$ of $A$ will be called a {\em standard orthonormal basis}
for $A$ if $X$ is an orthonormal set with respect to this Euclidean
structure, if the linear span of $X$ is a dense $^*$-subalgebra of $A$
(with respect to the $C^*$-norm), and if $1\in X$. The set difference
$X\backslash\{1\}$ will often be denoted by $X^\circ$.

\vspace{.5cm}

\noindent{\bf Lemma}  \ 
{\em Assume $A$ is a separable $C^*$-algebra and that $F\subseteq A$ is a
finite
orthonormal set containing 1. Then there exists a (countable) standard
orthonormal basis for $A$ which contains $F$.}

\vspace{.5cm}

\bproof
Choose a dense subset $\{a_1,a_2,a_3,\dots\}$ of $A$. Set
$X_0=F$. Construct inductively finite orthonormal sets
$X_0\subseteq X_1\subseteq X_2\subseteq \cdots$,
satisfying
\begin{itemize}
\item[(i)] $a_n\in \Span X_n$,
\item[(ii)] $\Span X_n$ is self-adjoint,
\item[(iii)] $x,y\in X_{n-1} \Rightarrow xy\in \Span X_n$,
\end{itemize}
for all $n\ge 1$, as follows. Suppose $X_{n-1}$ has been constructed. Let
$V_n$ be the finite dimensional subspace of $A$ spanned by $a_n$,
$X_{n-1}$, $X_{n-1}\cdot X_{n-1}$ and the adjoints of those
elements. Then choose $X_n$ to be an orthonormal basis for $V_n$ that
extends $X_{n-1}$.

It is now easily verified that $X=\bigcup_{n=0}^\infty X_n$ is a
standard orthonormal basis for $A$.
\eproof

\subsection{Reduced free products}
\label{subsec-2.2}
Let $(A_i,\tau_i)$, $i\in I$, be a family of unital $C^*$-algebras $A_i$
with faithful normalized traces $\tau_i$. To each such family one can
associate the reduced free product $C^*$-algebra
\[
(\frA,\tau)=\underset{i\in I}{*} (A_i,\tau_i),
\]
where $\frA$ is a unital $C^*$-algebra and $\tau$ is a
normalized faithful trace on $\frA$ (\cite{Voiculescu}, see also \cite{VKN}).

By construction, $A_i$ is a sub-$C^*$-algebra of $\frA$, and $\tau$
extends $\tau_i$ for each $i\in I$. Elements in $\frA$ of the form
\[
w = a_1a_2a_3\cdots a_n,
\]
where $a_j\in A_{i(j)}$, $\tau(a_j)=0$, and $i(1)\ne i(2),i(2)\ne
i(3),\dots,i(n-1)\ne i(n)$, are said to be {\em reduced words} of
(block-) length $n$, and $a_1,a_2,\dots,a_n$ are said to be the {\em
letters} of the word $w$. (It turns out that the block length is well
defined.) The unit 1 in $\frA$ is said to be a reduced word of length
0. For each reduced word $w$ of length $n\ge 1$ we have
$\tau(w)=0$.
The linear span of all reduced words in $\frA$ is a norm dense
$^*$-subalgebra of $\frA$.

Suppose now that $X_i$ is a standard orthonormal basis for $A_i$. For
each $n\ge 1$, let $Y_n$ be the set of all reduced words $x_1x_2x_3\cdots
x_n$ where $x_j\in X^\circ_{i(j)}$, and $i(1)\ne i(2),i(2)\ne
i(3),\dots,i(n-1)\ne i(n)$. Set $Y_0=\{1\}$ and set
\[
\underset{i\in I}{*} X_i=\bigcup_{n=0}^\infty Y_n.
\]

>From the construction of the reduced free products (see
\cite{Voiculescu}) it is easily seen that $^*_{i\in I} X_i$ is an
orthonormal set (with respect to the Euclidean structure on $\frA$
arising from $\tau$). The linear span of $^*_{i\in I}X_i$ is a
$^*$-algebra (because each $\Span X_i$ is a $^*$-algebra). The closure of
the linear span of $^*_{i\in I} X_i$ contains all reduced words in
$\frA$, and is therefore equal to $\frA$. This shows that $^*_{i\in I}
X_i$ is a standard orthonormal basis for $\frA$.

\section{The main result}
\label{sec-3}

As in Section \ref{sec-2}, let $(A_i,\tau_i)$, $i\in I$, be a family of
unital $C^*$-algebras $A_i$ with faithful normalized traces $\tau_i$,
and with standard orthonormal bases $X_i\subseteq A_i$. Let
\[
(\frA,\tau) = \underset{i\in I}{*} (A_i,\tau_i)
\]
be the reduced free product $C^*$-algebra, and let $Y=*_{i\in I} X_i$ be
the standard orthonormal basis for $\frA$ defined in
\ref{subsec-2.2}. Let
\[
E_n: \Span Y\to \Span Y_n
\]
be the orthogonal projection.

We shall in the first lemma of this section describe the element
$E_n(vw)$, where $v\in Y_k$ and $w\in Y_\ell$ for some $k$, $\ell$ and
$n$. As in Section \ref{sec-2} we shall equip $\frA$ with the Euclidean
structure
\[
\langle a,b\rangle = \tau(b^*a),\qquad \|a\|_2=\langle
a,a\rangle^{\frac12},\qquad a,b\in\frA.
\]

\begin{lemma}
\label{lemma-3.1}
Let $v\in Y_k$, let $w\in Y_\ell$ and let $n\ge 0$ be given.
\begin{itemize}
\item[(i)] Assume $|k-\ell|<n\le k+\ell$. Let $q$ be the integer
satisfying $k+\ell-n=2q$ or $k+\ell-n=2q+1$ (which entails that $0\le
q<\min\{k,\ell\}$). Write
\begin{alignat*}{4}
v &= v_1 x v_2, & \qquad v_1 &\in Y_{k-q-1}, &\qquad x &\in X^\circ_i,
&\qquad v_2 &\in Y_q,\\
w &= w_2yw_1, &\qquad w_1 &\in Y_{\ell-q-1}, &\qquad y &\in X^\circ_j, &\qquad
w_2 &\in Y_q.
\end{alignat*}
It follows that
\[
E_n(vw) = \left\{\begin{array}{cl} \langle v_2w_2,1\rangle v_1
xyw_1, &\qquad  \mbox{if $i\ne j$}\\
0, & \qquad  \mbox{if $i=j$} \end{array}\right.
\]
if $k+\ell-n$ is even, and
\[
E_n(vw) = \left\{\begin{array}{cl} \sum_{u\in X^\circ_i} \langle
v_2w_2,1\rangle \langle xy,u\rangle v_1uw_1, & \qquad  \mbox{if $i=j$}\\
0, & \qquad \mbox{if $i\ne j$} \end{array}\right.
\]
if $k+\ell-n$ is odd. (Observe that $\langle xy,u\rangle \ne 0$ for at
most finitely many $u\in X_i^\circ$ because $xy\in \Span X_i$.)
\item[(ii)] Assume $n=|k-\ell|$. Put $q=\min\{k,\ell\}$, so that
$k+\ell-n=2q$, and write
\begin{alignat*}{3}
v &= v_1v_2, &\qquad v_1 &\in Y_{k-q}, &\qquad v_2 & \in Y_q,\\
w &= w_2w_1, &\qquad w_1 &\in Y_{\ell-q}, &\qquad w_2 &\in Y_q.
\end{alignat*}
It follows that $v_1=1$ or $w_1=1$, and $E_n(vw)=\langle
v_2w_2,1\rangle v_1w_1$.
\item[(iii)] If $n<|k-\ell|$ or if $n>k+\ell$, then $E_n(vw)=0$.
\end{itemize}
\end{lemma} 

\bproof
We prove (i), (ii) and (iii) simultaneously by induction on
$\min\{k,\ell\}$. If $\min\{k,\ell\}=0$, then $v=1$ or $w=1$, and either
$n=|k-\ell|$ and $q=0$, or $n<|k-\ell|$, or $n>k+\ell$. The claims are
trivial in all six cases.

Consider now the case where  $\min\{k,\ell\}\ge 1$. Write $v=v'x'$ and
$w=y'w'$ with $v'\in Y_{k-1}$, $w'\in Y_{\ell-1}$, $x'\in X_s^\circ$ and
$y'\in X_t^\circ$. If $s\ne t$, then $vw$ is reduced, and so
\[
E_n(vw) = \left\{\begin{array}{cl} vw, &\qquad \mbox{if $n=k+\ell$}\\0,
&\qquad \mbox{if $n\ne k+\ell$.} \end{array}\right.
\]
This formula agrees with (iii). If $n=k+\ell$, then $q=0$ in (i) and
(ii), which entails $v_2=w_2=1$ and thereby $\langle v_2w_2,1\rangle
=1$. If $|k-\ell|\le n <k+\ell$, then $q\ge 1$ in (i) and (ii), and
$v_2w_2$ is a reduced word (because $s\ne t$). Hence $\langle
v_2w_2,1\rangle =0$. In either event, the expression for $E_n(vw)$
displayed above agrees with the formulae in (i) and (ii).

Suppose now that $s=t$. Then
\[
vw = \langle x'y',1\rangle v'w' + \sum_{u\in X_s^\circ} \langle
x'y',u\rangle v'uw'.
\]
Hence $E_n(vw)$ is as claimed in the lemma when $n\ge k+\ell-1$.

Consider now the case where $|k-\ell|\le n<k+\ell-1$. Then $q\ge 1$, and
in the notation of (i) and (ii) we can write $v_2=v'_2x'$ and
$w_2=y'w'_2$ for some $v'_2,w'_2\in Y_{q-1}$. Hence $v'=v_1xv'_2$ and
$w'=w'_2yw_1$. Now,
\[
E_n(vw) = \langle x'y',1\rangle E_n(v'w'),
\]
and $E_n(v'w')$ is by the induction hypothesis given by the formulae in
(i) and (ii). Since
\begin{eqnarray*}
\langle v_2w_2,1\rangle &=& \tau(v'_2x'y'w'_2)\\
&=& \tau\big(v'_2\langle x'y',1\rangle 1\cdot w'_2\big) +
\tau\big(v'_2(x'y'-\langle
x'y',1\rangle 1)w'_2\big)\\
&=& \tau\big(v'_2\langle x'y',1\rangle 1\cdot w'_2\big)\\
&=& \langle x'y',1\rangle \langle v'_2w'_2,1\rangle,
\end{eqnarray*}
the formulae for $E_n(vw)$ in (i) and (ii) are verified.

Finally, if $n<|k-\ell|$, then $n<|(k-1)-(\ell-1)|$, whence
\[
E_n(vw) = \langle x'y',1\rangle E_n(v'w')=0.
\] 
\eproof

\begin{definition}
\label{def-3.2}
{\rm 
For every $a\in\Span Y$ and for every $i\in I$ define $F_i(a)$ to be the
set of all $x\in X_i^\circ$ that appear as letters in words $w\in Y$ in the
support of $a$. (An element $w\in Y$ is said to lie in the support of $a$ if
$\langle a,w\rangle \ne 0$.)

Since the support of $a$ is finite, it follows that each $F_i(a)$ is
finite and that $F_i(a)\ne\emptyset$ only for finitely many $i\in I$.

Set
\[
K(a) = \max_{i\in I} \big(\sum_{x\in F_i(a)}\|x\|^2\big)^{\frac{1}{2}}
\]}
\end{definition}

\begin{lemma}
\label{lemma-3.3}
Let $a\in\Span Y_k$, let $b\in \Span Y_\ell$, and let $n\ge 0$ be
given. If $|k-\ell|\le n\le k+\ell$, then

\[
\|E_n(ab)\|_2 \le \left\{\begin{array}{rcl} \|a\|_2\; \|b\|_2&, &\qquad
\mbox{if $k+\ell-n$ is even}\\
K(a)\|a\|_2\; \|b\|_2&, &\qquad \mbox{if $k+\ell-n$ is odd.}
\end{array}\right.
\]

\noindent If $n<|k-\ell|$ or if $n>k+\ell$, then $E_n(ab)=0$.
\end{lemma}

\bproof
It is an immediate consequence of Lemma \ref{lemma-3.1} (iii) that
$E_n(ab)=0$ when $n<|k-\ell|$ and when $n>k+\ell$. Assume therefore that
$|k-\ell|\le n\le k+\ell$. Consider first the case where $k+\ell -n$ is
even. Write $k+\ell -n=2q$, and note that $0\le q\le
\min\{k,\ell\}$. Write
\[
a = \sum_{v_1,v_2} \alpha_{v_1v_2}v_1v_2, \qquad b=\sum_{w_1,w_2}
\beta_{w_2w_1}w_2w_1
\]
summing over all $v_1\in Y_{k-q}$, $v_2\in Y_q$ such that $v_1v_2$ is
reduced, respectively, over all $w_1\in Y_{\ell-q}$, $w_2\in Y_q$ such
that $w_2w_1$ is reduced. (Only finitely many $\alpha_{v_1v_2}$ and
$\beta_{w_2w_1}$ are non-zero.)

By Lemma  \ref{lemma-3.1},
\[
E_n(ab) = \sum_{v_1,w_1} \sum_{v_2,w_2} \alpha_{v_1v_2}\beta_{w_2w_1}
\langle v_2w_2,1\rangle v_1w_1
\]
summing over all $v_1\in Y_{k-q}$, $w_1\in Y_{\ell-q}$ and $v_2,w_2\in
Y_q$ such that the words $v_1v_2,w_2w_1$ and $v_1w_1$ are all
reduced. Hence
\[
\|E_n(ab)\|_2^2 \le \sum_{v_1,w_1} \big| \sum_{v_2,w_2}
\alpha_{v_1v_2}\beta_{w_2w_1} \langle v_2w_2,1\rangle \big|^2,
\]
now summing over all $v_1\in Y_{k-q}$, $w_1\in Y_{\ell-q}$ and
$v_2,w_2\in Y_q$ such that $v_1v_2$ and $w_2w_1$ are reduced. We use the
Cauchy-Schwartz inequality to estimate the right-hand side:
\begin{eqnarray*}
\big|\sum_{v_2,w_2}\alpha_{v_1v_2}\beta_{w_2w_1}\langle v_2w_2,1\rangle
\big|^2 &=& \big|\langle \sum_{w_2}
\beta_{w_2w_1}w_2,\sum_{v_2}\bar{\alpha}_{v_1v_2}v^*_2\rangle\big|^2\\
& \le &
\big\|\sum_{w_2}\beta_{w_2w_1}w_2\big\|^2_2\cdot\big\|
\sum_{v_2}\bar{\alpha}_{v_1v_2}v_2^*\big\|^2_2\\
&=& \sum_{w_2} |\beta_{w_2w_1}|^2\cdot \sum_{v_2}|\alpha_{v_1v_2}|^2.
\end{eqnarray*}
Hence
\begin{eqnarray*}
\|E_n(ab)\|^2_2 &\le &
\sum_{v_1,w_1}\sum_{w_2}|\beta_{w_2w_1}|^2\cdot\sum_{v_2}
|\alpha_{v_1v_2}|^2\\
&=& \sum_{v_1,v_2} |\alpha_{v_1v_2}|^2\cdot
\sum_{w_1,w_2}|\beta_{w_2w_1}|^2 = \|a\|^2_2\cdot \|b\|^2_2.
\end{eqnarray*}

Suppose now that $k+\ell -n$ is odd, and write $k+\ell -n=2q+1$ for some
integer $q$ satisfying $0\le q<\min\{k,\ell\}$. Write
\[
a = \sum_{i\in I}\sum_{v_1,x,v_2} \alpha_{v_1xv_2}v_1xv_2,\qquad
b=\sum_{i\in I}\sum_{w_1,y,w_2} \beta_{w_2yw_1}w_2yw_1
\]
summing over all $v_1\in Y_{k-q-1}$, $x\in X^\circ_i$, $v_2\in Y_q$ such
that $v_1xv_2$ is reduced, respectively, over all $w_1\in Y_{\ell-q-1}$,
$y\in X_i^\circ$, $w_2\in Y_q$ such that $w_2yw_1$ is reduced. By Lemma
\ref{lemma-3.1},
\[
E_n(ab) = \sum_{v_1,w_1}\sum_{i\in I}
\sum_{u\in X_i^\circ}\sum_{x,y\in X_i^\circ}
\sum_{v_2,w_2}
 \alpha_{v_1xv_2}\beta_{w_2yw_1}
\langle v_2w_2,1\rangle \langle xy,u\rangle v_1uw_1,
\]
summing over all $v_1\in Y_{k-q-1}$, $w_1\in Y_{\ell-q-1}$, $v_2,w_2\in
Y_q$ such that $v_1xv_2$ and $w_2yw_1$ are reduced. Hence
\[
\|E_n(ab)\|^2_2 = \sum_{v_1,w_1}\sum_{i\in I}\sum_{u\in X_i^\circ}
\big|
\sum_{x,y\in X_i^\circ}
\sum_{v_2,w_2}
\alpha_{v_1xv_2}\beta_{w_2yw_1}\langle v_2w_2,1\rangle \langle
xy,u\rangle\big|^2.
\]
For fixed $v_1,w_1$ and $i\in I$, put
\[
z = \sum_{x,y\in X_i^\circ} \big\langle \sum_{w_2} \beta_{w_2yw_1}w_2,
\sum_{v_2}\bar{\alpha}_{v_1xv_2}v_2^*\big\rangle xy \in \Span X_i.
\]
Then
\[
\big|
\sum_{x,y\in X_i^\circ}
\sum_{v_2,w_2}
\alpha_{v_1xv_2}\beta_{w_2yw_1}
\langle v_2w_2,1\rangle \langle xy,u\rangle \big|^2 = |\langle z,u
\rangle |^2,
\]
and since $\alpha_{v_1xv_2}=0$ if $x\notin F_i(a)$,
\begin{eqnarray*}
\|z\|^2_2 &=& \big\|\sum_{x\in F_i(a)} x \sum_{y\in X_i^\circ}
\big\langle \sum_{w_2}
\beta_{w_2yw_1}w_2,\sum_{v_2}\bar{\alpha}_{v_1xv_2}v_2^*\big\rangle
y\big\|_2^2\\
&\le & \big(\sum_{x\in F_i(a)} \|x\|\cdot \big\|\sum_{y\in
X_i^\circ}\big\langle \sum_{w_2}
\beta_{w_2yw_1}w_2,\sum_{v_2}\bar{\alpha}_{v_1xv_2}v_2^*\big\rangle
y\big\|_2\big)^2\\
&\le & \big(\sum_{x\in F_i(a)}\|x\|^2\big)\cdot\big(\sum_{x\in F_i(a)}
\big\|\sum_{y\in X_i^\circ}\big\langle \sum_{w_2}
\beta_{w_2yw_1}w_2,\sum_{v_2}\bar{\alpha}_{v_1xv_2}v_2^*\big\rangle
y\big\|^2_2\big)\\
& \le & K(a)^2  \sum_{x\in F_i(a)}\sum_{y\in X_i^\circ} \big|\big\langle
\sum_{w_2}\beta_{w_2yw_1}w_2,\sum_{v_2}\bar{\alpha}_{v_1xv_2}v_2^*
\big\rangle\big|^2\\
&\le & K(a)^2 \sum_{x,y\in X_i^\circ}
\big\|\sum_{w_2}\beta_{w_2yw_1}w_2\big\|^2_2\cdot
\big\|\sum_{v_2}\bar{\alpha}_{v_1xv_2}v_2^*\big\|^2_2\\
&=& K(a)^2 \sum_{x,y\in X_i^\circ} \sum_{w_2}
|\beta_{w_2yw_1}|^2\cdot\sum_{v_2} |\alpha_{v_1xv_2}|^2\\
&=& K(a)^2 \sum_{x,v_2} |\alpha_{v_1xv_2}|^2\cdot\sum_{y,w_2}
|\beta_{w_2yw_1}|^2.
\end{eqnarray*}
Hence
\begin{eqnarray*}
\sum_{u\in X_i^\circ} \big|
\sum_{x,y\in X_i^\circ}
\sum_{v_2,w_2}
\alpha_{v_1xv_2}\beta_{w_2yw_1}\langle v_2w_2,1\rangle
\langle xy,u\rangle\big|^2
&=& \sum_{u\in X_i^\circ} |\langle z,u\rangle |^2\le \|z\|_2^2\\
&\le & K(a)^2 \sum_{x,v_2}|\alpha_{v_1xv_2}|^2\cdot
\sum_{y,w_2}|\beta_{w_2yw_1}|^2.
\end{eqnarray*}
Finally, this proves that
\begin{eqnarray*}
\|E_n(ab)\|^2_2 &\le & \sum_{v_1,w_1} \sum_{i\in I} K(a)^2 \sum_{x,v_2}
|\alpha_{v_1xv_2}|^2\cdot\sum_{y,w_2}|\beta_{w_2yw_1}|^2\\
&=& K(a)^2 \big(\sum_{i\in I}\sum_{v_1,x,v_2}
|\alpha_{v_1xv_2}|^2\big)\big(\sum_{i\in I}\sum_{w_1,y,w_2}
|\beta_{w_2yw_1}|^2\big)\\
&=& K(a)^2 \|a\|^2_2 \|b\|^2_2.
\end{eqnarray*}
\eproof

\begin{lemma}
\label{lemma-3.4}
For each $a\in \Span Y_k$,
\[
\|a\|\le (2k+1)K(a)\|a\|_2.
\]
\end{lemma}

\bproof
It suffices to show that
\[
\|ab\|_2 \le (2k+1)K(a)\|a\|_2 \|b\|_2
\]
for all $b\in\Span Y$. Put $b_j=E_j(b)$. Then the following estimate
holds for each $n\ge 0$ by Lemma \ref{lemma-3.3}:
\begin{eqnarray*}
\|E_n(ab)\|_2 &=& \big\|\sum_{j=|n-k|}^{n+k} E_n(ab_j)\big\|_2\\
&\le & \sum_{j=|n-k|}^{n+k} \|E_n(ab_j)\|_2\\
&\le & \sum_{j=|n-k|}^{n+k} K(a)\|a\|_2 \|b_j\|_2\\
&\le & K(a) \|a\|_2(2k+1)^{\frac12} \big(\sum_{j=|n-k|}^{n+k}
\|b_j\|_2^2\big)^{\frac12}.
\end{eqnarray*}
Hence
\begin{eqnarray*}
\|ab\|^2_2 &\le & \sum^\infty_{n=0} \|E_n(ab)\|_2^2\\
&\le & (2k+1)K(a)^2\|a\|^2_2 \sum^\infty_{n=0} \sum_{j=|n-k|}^{n+k}
\|b_j\|^2_2\\
&\le & (2k+1)^2 K(a)^2 \|a\|^2_2 \sum^\infty_{j=0} \|b_j\|^2_2\\
&=& (2k+1)^2 K(a)^2 \|a\|^2_2 \|b\|^2_2.
\end{eqnarray*}
\eproof

\begin{lemma}
\label{lemma-3.5}
For each $a\in\Span\big(\bigcup^k_{j=0} Y_j\big)$,
\[
\|a\|\le (2k+1)^{\frac32} K(a) \|a\|_2.
\]
\end{lemma}

\bproof
Put $a_j=E_j(a)$. Observe that $K(a_j)\le K(a)$ (see Definition
\ref{def-3.2}). Lemma \ref{lemma-3.4} now yields
\begin{eqnarray*}
\|a\| &=& \big\|\sum^k_{j=0} a_j\big\|\le \sum^k_{j=0} \|a_j\|\\
&\le & \sum^k_{j=0} (2j+1)K(a_j)\|a_j\|_2\\
&\le & (2k+1)K(a)\sum^k_{j=0} \|a_j\|_2\\
&\le & (2k+1)K(a)(k+1)^{\frac12} \big(\sum^k_{j=0}
\|a_j\|^2_2\big)^{\frac12}\\
&=& (2k+1)(k+1)^{\frac12} K(a)\|a\|_2\\
&\le & (2k+1)^{\frac32} K(a)\|a\|_2.
\end{eqnarray*}
\eproof

\begin{lemma}
\label{lemma-3.6}
Suppose
\[
v = a_2a_2\cdots a_r,\quad w=b_1b_2\cdots b_s,\quad z=c_tc_{t-1}\cdots
c_1
\]
are three reduced words in $\frA$ (of block-length $r,s$ and $t$), and
suppose that $s<\min\{r,t\}$. Then $vwz$ is a linear combination of {\em
reduced} words of the form
\[
a_1a_2\cdots a_{r'}b'_1b'_2\cdots b'_s c_{t'}c_{t'-1}\cdots c_1
\]
and of (possibly unreduced) words of the form
\[
a_1a_2\cdots a_{r'}c_{t'}c_{t'-1}\cdots c_1,
\]
where $r'\ge r-s$ and $t'\ge t-s$ (in both cases).
\end{lemma}

\bproof
The proof is by induction on $s$. If $s=0$, then $w=1$ and
$vwz=a_1a_2\cdots a_rc_tc_{t-1}\cdots c_1$ in agreement with the
lemma. Let now $s>0$. Then
\[
a_r\in A_i,\quad b_1\in A_j,\quad b_s\in A_k,\quad c_t\in A_l
\]
for some $i,j,k,l\in I$. Consider the following four possibilities:

(i) $i\ne j$ and $k\ne l$, (ii) $i=j$ and $k\ne l$, (iii) $i\ne j$ and
$k=l$, and (iv) $i=j$ and $k=l$.

In case (i) the word $vwz$ is itself reduced and hence of the right
form. In case (iv), if $s\ge 2$, then set
\begin{alignat*}{2}
b'_1 &= a_rb_1 - \langle a_rb_1,1\rangle 1 \in A_i, & \qquad
\tau(b'_1) & =0,\\
b'_s &= b_sc_t - \langle b_sc_t,1\rangle 1\in A_k, & \qquad \tau(b'_s) & =0.
\end{alignat*}
If $s=1$, then $i=j=k=l$, and we set
\[
b'_1 = a_rb_1c_t-\langle a_rb_1c_t,1\rangle 1\in A_i,\quad \tau(b'_1)=0.
\]
If $s\ge 2$, then
\begin{eqnarray*}
vwz &=& a_1a_2\cdots a_{r-1}b'_1b_2\cdots b_{s-1}b'_sc_{t-1}\cdots
c_2c_1\\
&& + \langle a_rb_1,1\rangle a_1a_2\cdots a_{r-1}b_2\cdots
b_{s-1}b'_sc_{t-1}\cdots c_2c_1\\
&& + \langle b_sc_t,1\rangle a_1a_2\cdots a_{r-1}b'_1b_2\cdots
b_{s-1}c_{t-1}\cdots c_2c_1\\
&& + \langle a_rb_1,1\rangle \langle b_sc_t,1\rangle a_1a_2\cdots
a_{r-1}b_2\cdots b_{s-1}c_{t-1}\cdots c_2c_1,
\end{eqnarray*}
and if $s=1$, then
\[
vwz = a_1a_2\cdots a_{r-1}b'_1c_{t-1}\cdots c_2c_1 + \langle
a_rb_1c_t,1\rangle a_1a_2\cdots a_{r-1}c_{t-1}\cdots c_2c_1.
\]
The first term of each of these two expressions is reduced, and the
remaining three terms of the first expression are, by the induction
hypothesis, linear combinations of words of the desired form.

Cases (ii) and (iii) can be treated in a similar way.
\eproof

\begin{lemma}
\label{lemma-3.7}
Assume that for some distinct pair of indices $i_1,i_2\in I$ there
exist at least one unitary element in $X^\circ_{i_1}$ and at least two
unitary elements in $X^\circ_{i_2}$. Then for each $a\in \Span Y$ there
exist unitaries $u,v\in \Span Y$ and a constant $K<\infty$ such that
\[
\|(uav)^n\|_2 = \|a\|_2^n,\qquad
K\big((uav)^n\big) \le  K,
\]
for all $n\ge 1$.
\end{lemma}

\bproof
Let $x\in X^\circ_{i_1}$ and $y,z\in X^\circ_{i_2}$ be distinct unitary
elements. Let $k$ be the length of the longest word $w\in Y$ in the
support of $a$, so that $a\in \Span (\bigcup^k_{j=0} Y_j)$. Choose an
integer $l$ such that $l\ge (k+3)/2$, and set
\[
u' = (xy^*)^l,\qquad v=(xz)^l.
\]
Notice that $u',v\in\Span Y$ because $Y$ is a standard orthonormal basis.

We show that whenever $w\in Y_j$ and $j\le k$, then $u'wv$ is a linear
combination of reduced words in $Y$ starting with $x$ and ending with
$z$.

Let $u'_s$ and $v_r$ be the words consisting of the first $s$ letters of
$u'$, respectively, the last $r$ letters of $v$. It follows from Lemma
\ref{lemma-3.6} that $u'wv$ is a linear combination of reduced words of
the form $u'_sw'v_r$ and of possibly unreduced words of the form $u'_sv_r$
where $s,r\ge 2l-j\ge 3$ in both cases. Moreover, by the proof of Lemma
\ref{lemma-3.6}, and since $w\in Y$, each $w'$ above will belong to
$\Span Y$, whence also $u'_sw'v_r$ belongs to $\Span Y$. It follows that
$u'_sw'v_r$ is a linear combination of elements in $Y$ starting with $x$
and ending with $z$.

We must also show that $u'_sv_r$ is a linear combination of words in $Y$
starting with $x$ and ending with $z$, whenever $s,r\ge 3$. If $s,r$
either both are even or both are odd, then $u'_sv_r$ is reduced and
therefore expressible as a linear combination of the desired form. If
$s$ is even and $r$ is odd, then
\[
u'_sv_r = u'_{s-1}(y^*z)v_{r-1};
\]
and if $s$ is odd and $r$ is even, then
\begin{eqnarray*}
u'_sv_r &=& u'_{s-2} y^*xxzv_{s-2}\\
&=& u'_{s-2}y^*(x^2-\langle x^2,1\rangle 1)zv_{s-2}+\langle x^2,1\rangle
u'_{s-2}(y^*z)v_{s-2}.
\end{eqnarray*}
Hence $u'_sv_r$ are linear combinations of words (in $Y$) beginning with
$x$ and ending with $z$.

We have established that
\[
u'av = \sum^N_{j=1} \alpha_jw_j,
\]
where $w_1,w_2,\dots,w_N$ are distinct elements of $Y$ each starting
with $x$ and ending with $z$, and where each $w_j$ has length no greater than
$2l+k$.
Choose an integer $m$ such that $m\ge(2l+k+1)/2$ and consider the unitary
element of $Y$,
\[ r=(xy)(xz)^m(xy). \]
For each $n\ge1$ and each choice of $j_1,\ldots,j_n\in\{1,\ldots,N\}$, clearly
$rw_{j_1}rw_{j_2}\cdots rw_{j_n}$ is a reduced word and an element of $Y$.
Moreover, by the choice of $r$, if
\[ rw_{i_1}rw_{i_2}\cdots rw_{i_n}=rw_{j_1}rw_{j_2}\cdots rw_{j_n} \]
for some $n\ge1$ then $i_1=j_1$, $i_2=j_2$, $\ldots,$ $i_n=j_n$.

Let $u=ru'$.
We have shown that for each $n\ge 1$,
\[
(uav)^n = \sum_{j_1=1}^N \sum_{j_2=1}^N\cdots \sum_{j_n=1}^N
\alpha_{j_1}\alpha_{j_2}\cdots\alpha_{j_n}rw_{j_1}rw_{j_2}\cdots rw_{j_n},
\]
and the words $rw_{j_1}rw_{j_2}\cdots rw_{j_n}$ are reduced and distinct
elements of $Y$. The expression above is therefore the
(unique) way to write $(uav)^n$ as a linear combination of basis
elements in $Y$. We conclude that
\[
K\big((uav)^n\big) = K(uav)
\]
for all $n\ge 1$ (c.f.\ Definition \ref{def-3.2}), and so we may take
$K$ to be $K(uav)$. Also,
\begin{eqnarray*}
\|(uav)^n\|_2 &=& \sum^N_{j_1=1}\sum^N_{j_2=1}\cdots\sum^N_{j_n=1}
|\alpha_{j_1}\alpha_{j_2}\cdots\alpha_{j_n}|^2\\
&=& \sum^N_{j_1=1} |\alpha_{j_1}|^2 \cdot
\sum^N_{j_2=1}|\alpha_{j_2}|^2\cdot \cdots \cdot \sum^N_{j_n=1}
|\alpha_{j_n}|^2 = \|a\|_2^n.
\end{eqnarray*}
\eproof

Let $A$ be a unital $C^*$-algebra. Denote by $\U(A)$ and $\GL(A)$ the
group of unitary, respectively, invertible elements of $A$. For each
$a\in A$, let $r(a)$ denote the spectral radius of $A$.

If $u,v\in \U(A)$, then $r(uav)=r(vuavv^*)=r(vua)$. This shows that
\begin{equation}
\label{eq-3.1}
\inf_{u,v\in \U(A)} r(uav) = \inf_{u\in \U(A)} r(ua).
\end{equation}
Since
\begin{eqnarray*}
r(a) &\ge& \dist\big(a,\{a-\lambda 1\mid \lambda\in\bC\}\cap \GL(A)\big)\\
&\ge & \dist\big(a,\GL(A)\big),
\end{eqnarray*}
and since $\dist(ua,\GL(A))=\dist(a,\GL(A))$, we have
\begin{equation}
\label{eq-3.2}
\dist\big(a,\GL(A)\big) \le \inf_{u\in \U(A)}r(ua).
\end{equation}

If $\GL(A)$ is a dense subset of $A$ (with respect to the $C^*$-norm),
then $A$ is said to have {\em stable rank one}, written $\sr(A)=1$.

\begin{theorem}
\label{thm-3.8}
Let $(A_i,\tau_i)$ be a family of unital $C^*$-algebras $A_i$ with
faithful normalized traces $\tau_i$. Assume that for some distinct pair
of indices $i_1,i_2\in I$ there exist unitary elements $x\in A_{i_1}$,
$y,z\in A_{i_2}$ such that
\[
0 = \tau_{i_1}(x) = \tau_{i_2}(y) = \tau_{i_2}(z) = \tau_{i_2}(z^*y),
\]
(i.e.\ $\{1,x\}$ and $\{1,y,z\}$ are orthogonal sets when $A_{i_1}$ and
$A_{i_2}$ are equipped with the Euclidean structure arising from the
traces $\tau_{i_1}$ and $\tau_{i_2}$).

Let
\[
(\frA,\tau) = \underset{i\in I}{*} (A_i,\tau_i)
\]
be the reduced free product $C^*$-algebra.
Then $\frA$ has stable rank one.
\end{theorem}

The condition on the family $(A_i,\tau_i)$, $i\in I$, that there exist
unitaries $x,y$ and $z$ with the properties stated in the theorem was
considered by Avitzour in \cite{Av}. He proved that his condition implies
that $\frA$ is simple, and that $\frA$ has the Dixmier property.
Our proof of Theorem \ref{thm-3.8} does not rely on Avitzour's theorem.

The Avitzour condition will be investigated in more detail in Section
\ref{sec-4}, and we shall prove that it is implied by some rather
general conditions (see Proposition \ref{prop-4.1}).

\vspace{.5cm}

\noindent{\bf Proof of Theorem \ref{thm-3.8}:}
Each element
in $\frA$ can be approximated by elements belonging to $^*_{i\in I'}B_i$
for suitable separable
sub-$C^*$-algebras $B_i$ of $A_i$ (where $B_{i_1}$
can be assumed to contain $x$, and $B_{i_2}$ can be assumed to contain
$y,z$).
We may therefore assume that each $A_i$ is separable.

By Lemma \ref{subsec-2.1} each $A_i$ has a standard orthonormal basis
$X_i$ such that $x\in X^\circ_{i_1}$  and $y,z\in X^\circ_{i_2}$.
Let $Y$ denote the standard orthonormal basis $^*_{i\in I} X_i$ for $\frA$
(c.f.\ \ref{subsec-2.2}).
We will prove that
\begin{equation}
\label{eq-3.3}
\inf_{u\in \U(\frA)} r(ua) \le \|a\|_2 \quad (= \tau(a^*a)^{\frac12})
\end{equation}
for each $a\in\Span Y$.
For each such $a$ there exists a $k\in\bN$
such that
\[
a\in\Span\big(\bigcup^k_{j=0} Y_j\big).
\]
Find unitaries $u,v\in\Span Y$ and a constant $K<\infty$ as in Lemma
\ref{lemma-3.7}. Let $l\in\bN$ be large enough so that
\[
u,v\in\Span\big(\bigcup^l_{j=0} Y_j\big).
\]
Then, for each $n\ge 1$,
\[
(uav)^n \in\Span\big(\bigcup_{j=0}^{n(k+2l)} Y_j\big).
\]

Lemma \ref{lemma-3.5} and Lemma \ref{lemma-3.7} yield
\[
\|(uav)^n\| \le  K\big(2n(k+2l)+1\big)^{\frac32} \big\|(uav)^n\big\|_2 =
K(2n(k+2l)+1)^{\frac32} \|a\|_2^n.
\]
From (\ref{eq-3.1}) we get
\begin{eqnarray*}
\inf_{u\in \U(\frA)} r(ua) &\le & r(uav)\\
&=& \liminf_{n\to\infty} \;\big\|(uav)^n\big\|^{\frac{1}{n}}\\
&\le & \liminf_{n\to \infty}\: K^{\frac{1}{n}}
\big(2n(k+2l)+1\big)^{\frac{3}{2n}}\|a\|_2 = \|a\|_2.
\end{eqnarray*}

We now proceed to show that $\sr(\frA)=1$. If $\sr(\frA)\ne 1$, then by
\cite[Theorem 2.6]{Ro} there would exist an element $b\in\frA$ with
\[
1=\|b\|=\dist(b,\GL(\frA)).
\]
Now $b$ is the limit in norm of a sequence, $a_k$, of elements of $\Span Y$.
By (\ref{eq-3.2}) and (\ref{eq-3.3}) we have for each $k$, that
\[ \dist(a_k,\GL(A))\le\|a_k\|_2, \]
hence the same holds for $b$ and thus
$1=\|b\|=\|b\|_2$. Hence $0=\tau(1-b^*b)=\tau(1-bb^*)$, and $1-b^*b$ and
$1-bb^*$ are positive. Since $\tau$ is faithful, we conclude that $b$ is
unitary. But then $b\in\GL(\frA)$, contradicting $\dist(b,\GL(\frA))=1$.
\eproof

\begin{corollary}
\label{cor-3.9}
Let $G$ be a discrete group and suppose that $G$ is a free product
$G_1*G_2$ of two groups $G_1$ and $G_2$ satisfying $|G_1|\ge 2$
and $|G_2|\ge 3$. Then $C^*_{\red}(G)$ has stable rank one.
\end{corollary}

\bproof
It follows from \cite{Voiculescu} that $C^*_{\red}(G)$ is isomorphic to the
reduced free product
\[
\big(C^*_{\red}(G_1),\tau_1\big) * \big(C^*_{\red}(G_2),\tau_2\big),
\]
where $\tau_1$ and $\tau_2$ are the canonical traces. Since the group
elements in $G_j$ form an orthonormal set in $C^*_{\red}(G_j)$ with the
Euclidean structure arising from the trace $\tau_j$, we see that the
conditions of Theorem \ref{thm-3.8} are satisfied.
\eproof

Observe that Corollary \ref{cor-3.9} shows that $C^*_{\red}(F_n)$, where
$2\le n\le\infty$ and where $F_n$ is the free group on $n$ generators,
and the Choi algebra $C^*_{\red}(\bZ_2*\bZ_3)$ and its generalizations
$C^*_{\red}(\bZ_n*\bZ_m)$, where $n\ge 2$ and $m\ge 3$, all have stable
rank one.

\section{The Avitzour condition}
\label{sec-4}
\setcounter{equation}{0}

Given an integer $n\ge 2$. We shall derive some partial results
describing those unital $C^*$-algebras $A$, with a normalized trace
$\tau$, that contain unitary elements $u_1=1,u_2,\dots,u_n$ so that
$\{u_1,u_2,\dots,u_n\}$ is an orthonormal set with respect to $\tau$,
i.e.\ $\tau(u_j^*u_i)=0$ for $i\ne j$. In view of Theorem \ref{thm-3.8}
this will be of most interest for us when $n=2$ or $n=3$.

We state the positive results in the proposition below. Suppose $B$ is a
unital abelian sub-$C^*$-algebra of $A$ and let $\Bhat$ denote the
spectrum of $B$ (so that $B$ is $^*$-isomorphic to $C(\Bhat)$). The
restriction of the trace $\tau$ to $B$ is represented by a (regular)
Borel probability measure $\mu$ on $\Bhat$, and this measure appears in
part (i) of Proposition \ref{prop-4.1}.

When we call a sub-$C^*$-algebra of $A$ unital, then it is assumed to
contain the unit of $A$.

\begin{proposition}
\label{prop-4.1}
Let $A$ be a unital $C^*$-algebra with a normalized trace $\tau$, and
consider the Euclidean structure on $A$ defined by $\langle a,b\rangle
=\tau(b^*a)$, $a,b\in A$.

\begin{itemize}
\item[(i)] Suppose $A$ contains a unital abelian sub-$C^*$-algebra $B$
so that the measure $\mu$ on $\Bhat$ representing $\tau|B$ is diffuse
(i.e.\ has no atoms). Then $A$ contains a unitary element $u$ such that
$\{u^n\}_{n\in\bZ}$ is an orthonormal set in
$A$, i.e.\ $\tau(u^n)=0$ for each $n\ne 0$.
\item[(ii)] Suppose $A$ contains mutually orthogonal projections
$p_1,p_2,\dots,p_n$ with sum equal to 1, and suppose that
$\tau(p_1)=\tau(p_2)=\cdots =\tau(p_n)=1/n$. Then there is a unitary
element $u\in A$ such that $u^n=1$ and the set
$\{1,u,u^2,\dots,u^{n-1}\}$ is orthonormal in
$A$, i.e.\ $\tau(u^k)=0$ for $1\le k\le
n-1$.
\item[(iii)] Suppose $A$ contains mutually orthogonal projections
$p_1,p_2,\dots,p_n$ with sum equal to 1, and suppose that $\tau(p_j)\le
1/2$ for all $j$. Then there is a unitary $u$ in $A$ with $\tau(u)=0$,
and hence $\{1,u\}$ is an orthonormal set.
\item[(iv)] If $A$ contains a finite dimensional unital
sub-$C^*$-algebra that has no direct summand isomorphic to $\bC$, then
$A$ contains unitary elements $u,v$ such that the set $\{1,u,v\}$ is
orthonormal.
\end{itemize}
\end{proposition}

\bproof
Part (i) follows immediately from Lemma \ref{lemma-4.2} below. To prove
part (ii) put $\w=\exp(2\pi i/n)$ and set
\[
u = p_1+\w p_2+\cdots + \w^{n-1}p_n.
\]
Then $u$ has the desired properties. Under the assumptions of (iii) we
can find projections
\[
q_1,q_2,q_3 \in \Span\{p_1,p_2,\dots,p_n\}
\]
such that $q_1+q_2+q_3=1$ and $\tau(q_3)\le \tau(q_2)\le\tau(q_1)\le
1/2$. Observe that $\tau(q_2+q_3)\ge \tau(q_1)$, and that
\[
\{|\tau(q_1+\lambda q_2)|: \lambda\in\bT\} =
[\tau(q_1)-\tau(q_2),\tau(q_1)+\tau(q_2)].
\]
Hence there exist $\lambda,\mu\in\bT$ with $|\tau(q_1+\lambda
q_2)|=\tau(q_3)$ and $\tau(q_1+\lambda q_2)=-\mu \tau(q_3)$. Thus
$u=q_1+\lambda q_2+\mu q_3$ is a unitary element in $A$ with
$\tau(u)=0$.

(iv) It suffices to show that there for each $n\ge 2$ exist unitaries
$u,v\in M_n(\bC)$ such that $\{1,u,v\}$ is an orthonormal set with
respect to the normalized trace on $M_n(\bC)$. This follows from (ii)
when $n\ge 3$. For $n=2$ we can use
\[
u = \left(\begin{array}{cc} 1&0\\ 0&-1\end{array}\right), \qquad
v = \left(\begin{array}{cc} 0&1\\ 1&0\end{array}\right).
\]
\eproof

\begin{lemma}
\label{lemma-4.2}
Let $X$ be a compact Hausdorff space, and let $\mu$ be a diffuse Borel
probability measure on $X$. Then there exists a continuous function $u:
X\to\bT$ such that $E\mapsto\mu(u^{-1}(E))$, for $E$  a Borel subset of
$\bT$, is the Haar measure on $\bT$.
\end{lemma}

\bproof
For each $\e>0$ set
\begin{eqnarray*}
Y_\e &=& \{ f\in C(X,\bR) \mid \forall t\in\bR:
\mu(f^{-1}(\{t\}))<\e\},\\
Y &=& \bigcap_{n=1}^\infty Y_{1/n} = \{ f\in C(X,\bR)\mid \forall t\in\bR:
\mu(f^{-1}(\{t\}))=0\}.
\end{eqnarray*}
We begin by proving that each $Y_\e$ is an open and dense subset of
$C(X,\bR)$ --- with the uniform topology induced by $\|\cdot\|_\infty$
--- and consequently that $Y$ is a dense (and hence non-empty)
$G_\delta$-set.

Let $\e>0$. We show that each $f\in Y_\e$ is an inner point in
$Y_\e$. For each $t\in\bR$,
\begin{equation}
\label{eq-4.1}
\lim_{\delta\to 0^+} \mu\big(f^{-1}\big(\;]t-\delta,t+\delta[\;\big) =
 \mu(f^{-1})\big(\{t\}\big)<\e,
\end{equation}
and so there exists $\delta>0$ (depending on $t$) such that
$\mu\big(f^{-1}\big(\;]t-\delta,t+\delta[\;\big)\big)<\e$. By compactness of
$f(X)$ there
exist $t_1,t_2,\dots,t_n\in\bR$ and $\delta_1,\delta_2,\dots,\delta_n>0$
such that
\begin{gather*}
f(X)\subseteq \bigcup^n_{j=1}\; ]t_j-\delta_j,t_j+\delta_j[,\\
\mu(f^{-1}(]t_j-\delta_j,t_j+\delta_j[)) < \e,\qquad j=1,2,\dots,n.
\end{gather*}
There exists a $\delta>0$ with the property that for every $t\in f(X)$
the interval $]t-\delta,t+\delta[$ is contained in
$]t_j-\delta_j,t_j+\delta_j[$ for some $j$ (depending on $t$). It
follows that
\[
\mu(f^{-1}(]t-\delta,t+\delta[))<\e
\]
for every $t\in f(X)$, and hence for every $t\in\bR$.

Let $g\in C(X,\bR)$ and suppose that $\|f-g\|_\infty<\delta$. Then
\[
g^{-1}(\{t\}) \subseteq f^{-1}(]t-\delta,t+\delta[)
\]
for every $t\in\bR$. Hence $\mu(g^{-1}(\{t\}))<\e$ for every $t\in\bR$
and so $g\in Y_\e$. This shows that $f$ is an inner point of $Y_\e$.

We prove next that $Y_\e \subseteq \overline{Y_{\frac{2}{3}\e}}$ for
every $\e>0$. Since $C(X,\bR)=Y_2$, this implies that each $Y_\e$ is
dense in $C(X,\bR)$. Let $f\in Y_\e$ and let $r>0$. We shall find a
$g\in Y_{\frac{2}{3}\e}$ with $\|f-g\|_\infty <r$. Set
\[
\big\{t\in\bR \mid \mu(f^{-1}(\{t\}))\ge \frac23 \e \big\} =
\{t_1,t_2,\dots,t_n\}.
\]
(If this set is empty, then $f\in Y_{\frac{2}{3}\e}$ and we may take
$g=f$.)

Using (\ref{eq-4.1}) and that $\mu(f^{-1}(\{t\}))>0$ for at most
countably many $t\in\bR$ one can find a $\delta\in\;]0,\frac{r}{2}[$
satisfying
\begin{gather*}
\mu\big(f^{-1}(\;]t_j-\delta,t_j+\delta[\;)\big)<\e,\\
\mu\big(f^{-1}(\{t_j-\delta\})\big) = \mu\big(f^{-1}(\{t_j+\delta\})\big)=0,
\end{gather*}
for each $j=1,2,\dots,n$, and such that the closed intervals
$[t_1-\delta,t_1+\delta],\dots,[t_n-\delta,t_n+\delta]$ are mutually
disjoint.

Put
\[
K_j=f^{-1}([t_j-\delta,t_j+\delta]),\quad
L_j=f^{-1}(\{t_j-\delta\}),\quad M_j=f^{-1}(\{t_j+\delta\}).
\]
Since $\mu$ restricts to a diffuse finite Borel measure on $K_j$ there
exists a Borel subset $E_j$ of $K_j$ with $\mu(E_j)=\frac12
\mu(K_j)$. Put
\[
E'_j = (E_j\cup L_j)\backslash M_j,\quad E''_j=((K_j\backslash E_j)\cup
M_j)\backslash L_j.
\]
Then $\mu(E'_j)=\mu(E''_j)=\frac12\mu(K_j)$, $E'_j\cap E''_j=\emptyset$,
$L_j\subseteq E'_j$ and $M_j\subseteq E''_j$. Since $\mu$ is regular, we can
find, necessarily disjoint, compact sets $L'_j,M'_j$ satisfying
\begin{alignat*}{2}
L_j & \subseteq L'_j\subseteq E'_j,\qquad & \mu(L'_j) & \ge \frac13\mu(K_j),\\
M_j & \subseteq M'_j\subseteq E''_j,\qquad & \mu(M'_j) & \ge \frac13\mu(K_j).
\end{alignat*}
Find next continuous functions $g_j: K_j\to [t_j-\delta,t_j+\delta]$
with
\[
g_j|L'_j \equiv t_j-\delta,\quad g_j|M'_j \equiv t_j+\delta.
\]
Observe that $|f(x)-g_j(x)|<r$ for all $x\in K_j$ because $\delta<r/2$.

Define $g: X\to\bR$ to be
\[
g(x) = \left\{\begin{array}{ll} g_j(x), & \mbox{if $x\in K_j$,}\\
f(x), & \mbox{if $x\in X\backslash \bigcup^n_{j=1} K_j$.} \end{array}\right.
\]
Then $g\in C(X,\bR)$ and $\|f-g\|_\infty <r$. If $t\in
[t_j-\delta,t_j+\delta]$, then
\[
\mu(g^{-1}(\{t\})) = \mu(g_j^{-1}(\{t\}))\le \frac23 \mu(K_j) < \frac23
\e.
\]
If $t\in \bR\backslash \bigcup^n_{j=1}[t_j-\delta,t_j+\delta]$, then $t\notin
\{t_1,t_2,\dots,t_n\}$ and so
\[
\mu(g^{-1}(\{t\})) = \mu(f^{-1}(\{t\}))<\textstyle{\frac23} \e.
\]
Hence $g\in Y_{\frac{2}{3}\e}$ as desired.

It has now been proved that $Y$ is non-empty. Let $f$ be any function in
$Y$, and choose $a,b\in\bR$ such that $f(X)\subseteq [a,b]$. Define $g:
[a,b]\to [0,1]$ to be
\[
g(t) = \mu(f^{-1}([a,t])).
\]
Thus $g$ is continuous, surjective and increasing (although perhaps not
strictly increasing). Put $h=g\circ f: X\to [0,1]$, and observe that $h$
is continuous. Let $s\in [0,1]$ and set
\[
t = \max\{t'\in [a,b] \mid g(t')=s\}.
\]
Then $g^{-1}([0,s]) = [a,t]$, whence $h^{-1}([0,s])=f^{-1}([a,t])$ and
\[
\mu(h^{-1}([0,s])) = \mu(f^{-1}([a,t]))=g(t)=s.
\]
This shows that $\mu(h^{-1}(E))=m(E)$ for every Borel subset $E$ of
$[0,1]$, where $m$ is the Lebesgue measure. It follows that the function
\[
u(x) = e^{2\pi ih(x)},\quad x\in X,
\]
has the desired properties.
\eproof

We now turn to some negative results.

\begin{proposition}
\label{prop-4.3}
Suppose $X$ is a compact Hausdorff space and that $\mu$ is a Borel
probability measure on $X$ such that $\mu(\{x_0\})>\frac{1}{n}$ for some
$x_0\in X$. Then there do not exist continuous (or measurable)
functions $u_1,u_2,\dots,u_n: X\to\bT$ such that $\{u_1,u_2,\dots,u_n\}$
forms an orthonormal set in the Hilbert space $L^2(X,\mu)$.
\end{proposition}

\bproof
Assume to the contrary that $\mu(\{x_0\})=a> 1/n$ for some
$x_0\in X$ and that $\{u_1,u_2,\dots,u_n\}$ is an orthonormal set of
functions $u_1,u_2,\dots,u_n: X\to\bT$. Upon replacing each $u_j$ with
$\overline{u_j(x_0)}u_j$ we may assume that $u_1(x_0)=u_2(x_0)=\cdots =
u_n(x_0)=1$. Define a new inner product on $L^2(X,\mu)$ by
\[
\langle f,g\rangle = \frac{1}{1-a} \int_{X\backslash \{x_0\}}
f\bar{g} d\mu.
\]
Then
\[
\delta_{ij}=\int_X u_i\bar{u}_jd\mu = a+(1-a)\langle
u_i,u_j\rangle,
\]
which shows that
\[
\langle u_i,u_j\rangle = \left\{\begin{array}{ccl} 1 &, & \mbox{if $i=j$}\\
-\frac{a}{1-a} & , & \mbox{if $i\ne j$}. \end{array}\right.
\]

The proposition therefore follows from the linear algebra fact, proved
below, that if $v_1,v_2,\dots,v_n$ are unit vectors in an Euclidean
space, satisfying $\langle v_i,v_j\rangle =-\alpha$ when $i\ne j$, then
$\alpha\le (n-1)^{-1}$.

We prove the last claim by induction on $n$, the ground step $n=2$ being
trivial. Given unit vectors $v_1,v_2,\dots,v_n$ satisfying $\langle
v_i,v_j\rangle =-\alpha$ when $i\ne j$, set
\[
w_j = \frac{1}{\sqrt{1-\alpha^2}} \big(v_j-\langle v_j,v_n\rangle
v_n\big)=\frac{1}{\sqrt{1-\alpha^2}}(v_j+\alpha v_n),
\]
for $j=1,2,\dots,n-1$. Then $\|w_j\|=1$ and $\langle w_i,w_j\rangle
=-\alpha(1-\alpha)^{-1}$ when $i\ne j$. Hence $\alpha(1-\alpha)^{-1}\le
(n-2)^{-1}$ by the induction hypothesis, and this implies
$\alpha\le (n-1)^{-1}$.
\eproof

\begin{example}
\label{ex-4.4}
{\rm Let $m$ be the Lebesgue measure on $[0,1]$, and let $\delta$ be the
Dirac measure on $[0,1]$ supported on $\{0\}$. Set $\mu=\frac12\delta
+\frac12 m$ and define a (tracial) state $\tau$ on $C([0,1])$ by
\[
\tau(f) = \int_0^1 f(x)d\mu(x).
\]
Then $\mu(\{0\}) = \frac12$, but, as a routine calculation will show,
\[
\big\{\tau(u)\mid u\in \U(C([0,1]))\big\} =
\big\{\lambda\in\bC\mid 0<|\lambda|\le
1\big\}.
\]
In particular, there is no unitary element $u$ of $C([0,1])$ with $\tau(u)=0$.}
\end{example}

There is no (obvious) analogue of Proposition \ref{prop-4.1} (iii) that
works for $n=3$ as will be shown in Proposition \ref{prop-4.6} below.

\begin{lemma}
\label{lemma-4.5}
There is a positive real number $\alpha_0<\frac12$ so that the following
holds for each $\alpha\in [\alpha_0,\frac12]$. Let $\beta,\gamma\in\bR^+$ be
determined by the equations
\[
\alpha^2 - \beta\sqrt{1-\alpha^2} = -\alpha,\quad
\alpha^2+\beta^2+\gamma^2=1.
\]
Then each triple $\lambda_1,\lambda_2,\lambda_3\in\bC$ satisfying
\begin{gather}
\label{eq-4.2}
|-\alpha+\sqrt{1-\alpha^2}\lambda_j|=1, \\
\label{eq-4.3}
1-\sqrt{3}\gamma \le |-\alpha-\beta\lambda_j|\le 1+\sqrt{3}\gamma,
\end{gather}
for $j=1,2,3$, has $\lambda_1+\lambda_2+\lambda_3\ne 0$.
\end{lemma}

\bproof
If $\alpha=\frac12$, then $\beta=\sqrt{3}/2$ and $\gamma=0$. In this
case (\ref{eq-4.2}) and (\ref{eq-4.3}) become equivalent to
\[
|-\frac12 +
\frac{\sqrt{3}}{2}\lambda_j|=|-\frac12-\frac{\sqrt{3}}{2}\lambda_j|=1.
\]
These equations can be rewritten as
\[
|\lambda_j|^2 - \frac{2}{\sqrt{3}}
\Re(\lambda_j)=|\lambda_j|^2+\frac{2}{\sqrt{3}} \Re(\lambda_j)=1.
\]
Hence $|\lambda_j|=1$ and $\Re(\lambda_j)=0$, and so
$\lambda_1,\lambda_2,\lambda_3\in \{-i,i\}$.

There exist continuous functions $c_1,c_2: [0,\frac12]\to\bR^+$
satisfying $c_1(\frac12)=c_2(\frac12)=1$ and
\[
c_1(\alpha) \le |\Im(\lambda)|\le c_2(\alpha)
\]
for all $\lambda\in\bC$ satisfying (\ref{eq-4.2}) and
(\ref{eq-4.3}). (Concrete expressions for possible choices of such
functions $c_1$ and $c_2$ can be derived through a straightforward, but
rather tedious, calculation, which we shall omit.)

By the continuity of $c_1$ and $c_2$ there is a positive number
$\alpha_0<\frac12$ such that $c_1(\alpha)>\frac23$ and
$c_2(\alpha)<\frac43$ for all $\alpha\in [\alpha_0,\frac12]$. Now, if
$\lambda_1,\lambda_2,\lambda_3\in\bC$ satisfy (\ref{eq-4.2}) and
(\ref{eq-4.3}), then $\frac23 <|\Im (\lambda_j)|<\frac43$, $j=1,2,3$,
which entails $\lambda_1+\lambda_2+\lambda_3\ne 0$.
\eproof

The following proposition should be compared with (ii) and (iii) in
Proposition \ref{prop-4.1}.

\begin{proposition}
\label{prop-4.6}
There is a positive real number $\beta_0<\frac13$ so that the following
holds for each $\beta\in [\beta_0,\frac13]$. Let $X$ be the 4-point
space $\{1,2,3,4\}$ and let $\mu$ be the measure on $X$ given by
\[
\mu(\{1\}) = \mu(\{2\}) = \mu(\{3\}) = \frac13(1-\beta),\quad
\mu(\{4\})=\beta.
\]
Then there exists no pair of functions $u,v: X\to\bT$ satisfying
\begin{equation}
\label{eq-4.4}
\int_X ud\mu = \int_X vd\mu = \int_X u\bar{v}d\mu =0.
\end{equation}
\end{proposition}

\bproof
Let $\alpha_0$ be as in Lemma \ref{lemma-4.5}, and set
\[
\beta_0=\max\{\frac14,\alpha_0(1+\alpha_0)^{-1}\}.
\]
Then $\frac14\le\beta_0 <\frac13$. Let $\beta\in [\beta_0,\frac13]$ and
suppose, to reach a contradiction, that $u,v: X\to\bT$ satisfy
(\ref{eq-4.4}). Upon replacing $u$ and $v$ with $\overline{u(4)}u$ and
$\overline{v(4)}v$ we may assume that $u(4)=v(4)=1$.

Set $X_0=\{1,2,3\}$, let $u_1$ and $v_1$ be the restrictions of $u$ and
$v$ to $X_0$, and define an inner product on the 3-dimensional vector
space of all functions from $X_0$ to $\bC$ by
\[
\langle f,g\rangle = \frac13 \sum^3_{j=1} f(j)\overline{g(j)}.
\]
Then
\begin{eqnarray*}
0 &=& \int_X ud\mu = \beta+(1-\beta)\langle u_1,1\rangle,\\
0 &=& \int_X vd\mu = \beta + (1-\beta)\langle v_1,1\rangle,\\
0 &=& \int_X u\bar{v}d\mu = \beta + (1-\beta)\langle u_1,v_1\rangle.
\end{eqnarray*}
Set $\alpha=\beta(1-\beta)^{-1}$, and notice that
$\alpha\in[\alpha_0,\frac12]$. Moreover,
\[
\langle u_1,1\rangle = \langle v_1,1\rangle = \langle u_1,v_1\rangle
=-\alpha.
\]
Find $x,y: X_0\to\bC$ such that $\{1,x,y\}$ is an orthonormal set, and
\[
u_1 = -\alpha\cdot 1+\sqrt{1-\alpha^2} \cdot x,\qquad v_1 = -\alpha\cdot
1-\beta x+\gamma y,
\]
where $\beta,\gamma\in\bR^+$ are determined by
$\alpha^2-\beta\sqrt{1-\alpha^2}=-\alpha$ and
$\alpha^2+\beta^2+\gamma^2=1$. Since $|y(j)|^2 \le 3\langle y,y\rangle
=3$ for $j=1,2,3$, we conclude that
\begin{eqnarray*}
|-\alpha+\sqrt{1-\alpha^2}\cdot x(j)| &=& |u_1(j)|=1,\\
|-\alpha-\beta x(j)| &=& |v_1(j)-\gamma y(j)|\in
[1-\sqrt{3}\gamma,1+\sqrt{3}\gamma].
\end{eqnarray*}
Lemma \ref{lemma-4.5} now yields $x(1)+x(2)+x(3)\ne 0$ in contradiction
with
\[
0 = \langle x,1\rangle = \frac13 \sum^3_{j=1} x(j).
\]
\eproof

\section{Other structural results}
\label{sec-5}
\setcounter{equation}{0}

We shall in this section discuss structural results of reduced free
products
\[
(\frA,\varphi) = (A_1,\varphi_1)*(A_2,\varphi_2),
\]
where $A_1,A_2$ are unital $C^*$-algebras with faithful states
$\varphi_1$ and $\varphi_2$, in the cases that are not covered in
Theorem \ref{thm-3.8}. The discussion will mostly be recapitulations of
already known results. As we shall see, our knowledge is somewhat
fragmented, but it indicates the shape of a clear picture.

Let us first consider the question as to when $\frA$ is simple. Avitzour
has, as mentioned earlier, proved the following:

\begin{theorem} \label{thm-5.1} {\bf (Avitzour \cite[Proposition 3.1]{Av})}
Suppose there exist unitary elements $x\in A_1$ and $y,z\in A_2$ such
that $x$ belongs to the centralizer of $\varphi_1$, $y$ belongs to the
centralizer of $\varphi_2$, and
\[
0 = \varphi_1(x) = \varphi_2(y)=\varphi_2(z) = \varphi_2(z^*y).
\]
Then $\frA$ is simple, and for every $a\in\frA$, the element
$\varphi(a)1$ belongs to the closure of the convex hull of the set
$\{uau^*\mid u\in \U(\frA)\}$.
\end{theorem}

In the negative direction we have:

\begin{proposition}
\label{prop-5.2}
Suppose $p_1\in A_1$ and $p_2\in A_2$ are non-trivial central
projections such that $p_jA_jp_j=\bC p_j$, $j=1,2$, and
$\varphi_1(p_1)+\varphi_2(p_2)\ge 1$. Then $\frA$ is not simple.
\end{proposition}

\bproof
The irreducible representations of $C^*(1,p_1,p_2)$ are described in
\cite[Theorem 12]{ABH}. From that description (and upon replacing $p_1$
and $p_2$ by $1-p_1$ and $1-p_2$ if necessary) it follows that
$\|p_1p_2\|=1$ if (and only if) $\varphi_1(p_1)+\varphi_2(p_2)\ge
1$. Since
\[
A_j = \bC p_j +(1-p_j)A_j(1-p_j),
\]
it follows from \cite[Proposition 2]{ABH} that there exists a
multiplicative state $\rho: \frA\to\bC$ such that
$\rho(p_1)=\rho(p_2)=1$. Since $\frA\ne\bC$ (because $A_1\ne\bC$ or
$A_2\ne\bC$ by assumption), the kernel of $\rho$ is a non-trivial ideal
of $\frA$.
\eproof

\begin{remark}
\label{rem-5.4}
{\rm Suppose $(A_1,\varphi_1)$ and $(A_2,\varphi_2)$ satisfy the Avitzour
conditions (i.e.\ the conditions in Theorem \ref{thm-5.1}), so that
$\frA$ is simple.

If $\tau$ is a trace on $\frA$, then for each $a\in\frA$ and each $x$ in
the closure of the convex hull of $\{uau^*\mid u\in \U(\frA)\}$ we have
$\tau(x)=\tau(a)$. By Avitzour's theorem (\ref{thm-5.1}) this implies that
$\tau=\varphi$. Since $\varphi$ extends $\varphi_1$ and $\varphi_2$, we
conclude that $\varphi_1$ and $\varphi_2$ must both be traces. In other
words, if the Avitzour conditions are satisfied, then $\frA$ admits a
trace if and only if $\varphi_1$ and $\varphi_2$ are both traces. If
$\varphi_1$ and $\varphi_2$ are both traces, then $\frA$ has stable rank one by
Theorem \ref{thm-3.8}.

Suppose now that either $\varphi_1$ or $\varphi_2$ is not a trace. Then
$\frA$ does not admit a trace. To analyze this situation in more detail,
consider the following five degrees of being infinite, where $A$ is a
simple unital $C^*$-algebra.
\begin{enumerate}
\item $A$ admits no trace.
\item $A$ admits no quasi-trace.
\item $A$ is not stably finite.
\item $A$ is infinite.
\item $A$ is purely infinite.
\end{enumerate}
Some words of explanations: A projection in a $C^*$-algebra is called
infinite if it is Murray-von Neumann equivalent to a proper
subprojection of itself; and it is called finite otherwise. A unital
$C^*$-algebra is called finite, respectively infinite, if its unit, as a
projection, is finite, respectively infinite. (Equivalently, a
$C^*$-algebra is infinite if and only if it contains an infinite
projection.) If $M_n(A)$ is finite for all $n\ge 1$, then $A$ is said to
be stably finite. A $C^*$-algebra is purely infinite if all its non-zero
hereditary sub-$C^*$-algebras contain an infinite projection.

The following implications hold for all simple, unital $C^*$-algebras:
\[
1) \Leftarrow 2) \Leftrightarrow 3)\Leftarrow 4) \Leftarrow 5).
\]
The equivalence $2) \Leftrightarrow 3)$ is proved in \cite{BH}, and the
reader can find a definition of a quasi-trace the same place. All other
implications are trivial. It is proved in \cite{Ha1} that $1)\Rightarrow
2)$ whenever $A$ is exact. There is no known example of a simple unital
$C^*$-algebra without a trace, which is not purely infinite, i.e.\
there are no known counterexamples to the implication $1)\Rightarrow
5)$. If $4)\Rightarrow 5)$ were known to be true, then $3)\Rightarrow
4)$ would also be true.

If a $C^*$-algebra has a stable rank one, then it must be stably finite (by
\cite[Theorem 3.3]{Rf} and since each isometry has distance 1 to the
invertible elements). Hence if our reduced free product $\frA$ is exact
(in addition to our assumption that either $\varphi_1$ or $\varphi_2$ is
not a trace), then $\frA$ does not have stable rank one. Exactness might
be preserved under forming reduced free products. It is  known
that $C^*_{\red}(F_n)$ and $C^*_{\red}(\bZ_2*\bZ_3)$ are exact.

In certain concrete examples it is quite easy to see that a reduced free
product $C^*$-algebra is infinite. For example, if $A_1=C^*(\bZ)$ with
$\varphi_1$ the canonical trace and $A_2=M_2(\bC)$ with
\[
\varphi_2\left(\left(\begin{array}{cc} a_{11} & a_{12}\\a_{21} &
a_{22}\end{array}\right)\right) = \frac13 a_{11} + \frac23 a_{22},
\]
then $\frA=A_1*A_2$ is infinite. Indeed, let
$e= {1\,0\choose 0\,0}\in A_2$ and let $u\in
A_2$ be the unitary corresponding to a generator of $\bZ$. Then $e$ and
$f=u(1-e)u^*$ are free projections, and $\frac13 =
\varphi(e)<\varphi(f)=\frac23$. Hence $e$ is equivalent to a proper
subprojection of $f$ (cf.~\cite{ABH}). At the same time, $e$
and $f$ are equivalent, and this shows that $f$ is infinite.

It is in general much harder to show that a given infinite simple
$C^*$-algebra is purely infinite. 
Some examples of purely
infinite reduced free products (of finite $C^*$-algebras) are given in
\cite{DR}.}
\end{remark}


\begin{thebibliography}{99}

\bibitem{ABH} J.~Anderson, B.~Blackadar and U.~Haagerup, Minimal
projections in the reduced group $C^*$-algebras of $\bZ_n*\bZ_m$, {\em
J.\ Operator Theory} {\bf 26} (1991), 3--23.


\bibitem{Av} D.~Avitzour, Free products of $C^*$-algebras, {\em
Trans.\ Amer.\ Math.\ Soc.} {\bf 271} (1982), 423--435.

\bibitem{Bergman} G.M.~Bergman, Modules over coproducts of rings, {\em
Trans.\ Amer.\ Math.\ Soc.} {\bf 200} (1974), 1--32.

\bibitem{BH} B.~Blackadar and D.~Handelman, Dimension functions and
traces on $C^*$-algebras, {\em J.\ Funct.\ Anal.} {\bf 45} (1982), 297--340.

\bibitem{CohnZZFreeRings} P.M.~Cohn, {\em Free Rings and their Relations},
Academic Press, 1971.

\bibitem{DR} K.J.~Dykema and M.~R\o rdam, Purely infinite simple
$C^*$-algebras arising from free product constructions, to appear in {\em
Can.\ J.\ Math.}

\bibitem{FR} P.~Friis and M.~R\o rdam, Almost commuting self-adjoint
matrices --- a short proof of Huaxin Lin's theorem, to appear in {\em
J.\ reine angew.\ Math.}

\bibitem{Ha} U.~Haagerup, An example of a non nuclear $C^*$-algebra,
which has the metric approximation property, {\em Invent.\ Math.} {\bf
50} (1979), 279--293.

\bibitem{Ha1} U.~Haagerup, Quasi traces on exact $C^*$-algebras are
traces, in preparation.

\bibitem{PV} M.~Pimsner and D.~Voiculescu, K--groups of reduced crossed
products by free groups, {\em J.\ Operator Theory} {\bf 8} (1982), 131--156.

\bibitem{Rf} M.~Rieffel, Dimension and stable rank in the $K$-theory of
$C^*$-algebras, {\em Proc.\ London Math.\ Soc.} {\bf (3) 46} (1983), 301--333.

\bibitem{RieffelZZKthCrPr} M.~Rieffel, K--theory of crossed products of
C$^*$--algebras by discrete groups, {\em Group actions on rings, (Brunswick,
Maine, 1984)} Contemporary Math. 43, (1985), 253--265.

\bibitem{Rf1} M.~Rieffel, The homotopy groups of the unitary groups of
non-commutative tori, {\em J.\ Operator Theory} {\bf 17} (1987), 237--254.

\bibitem{Ro} M.~R\o rdam, Advances in the theory of unitary rank and
regular approximation, {\em Ann.\ of Math.} {\bf 128} (1988), 153--172.

\bibitem{Voiculescu} D.~Voiculescu, Symmetries of some reduced free
product $C^*$-algebras, {\em Operator Algebras and Their Connections with
Topology and Ergodic Theory}, Lecture Notes in Mathematics, Volume
1132, Springer-Verlag, 1985, 556--588.

\bibitem{VKN} D.\ Voiculescu, K.J.~Dykema, A.~Nica, {\em Free Random
Variables}, CRM Monograph Series vol.\ 1, American Mathematical Society, 1992.


\end{thebibliography}
\end{document}